\documentclass[10pt,a4paper]{amsart}
\usepackage{hyperref}
\usepackage{mathrsfs}
\usepackage{amsmath}
\usepackage{amssymb}
\usepackage{a4}
\usepackage[all]{xy}
\usepackage{latexsym}

\title{Embedding the flag representation in divided powers}
\author{Geoffrey M.L. Powell}
\address{Laboratoire Analyse, G\'eom\'etrie et Applications, UMR 7539\\ Institut Galil\'ee, Universit\'e Paris 13, 93430 Villetaneuse, France}
\email{powell@math.univ-paris13.fr} 
\date{}

\newtheorem{thm}{Theorem}[section]
\newtheorem{prop}[thm]{Proposition}

\newtheorem{lem}[thm]{Lemma}

\newtheorem{THM}{Theorem}

\newtheorem{PROP}[THM]{Proposition}

\theoremstyle{definition}
\newtheorem{defn}[thm]{Definition}

\theoremstyle{remark}
\newtheorem{rem}[thm]{Remark}
\newtheorem{nota}[thm]{Notation}

\long\def\forget#1{}

\newcommand{\fl}{{\Phi}}
\newcommand{\qpol}{\mathfrak{q}}
\newcommand{\flag}{\mathfrak{Flag}}

\renewcommand{\hom}{\mathrm{Hom}}

\newcommand{\zed}{\mathbb{Z}}

\newcommand{\field}{\mathbb{F}}
\newcommand{\f}{\mathscr{F}}

\newcommand{\op}{^\mathrm{op}}

\newcommand{\sym}{\mathfrak{S}}

\begin{document}
\maketitle

\begin{abstract}
A generalization of a theorem of Crabb and Hubbuck concerning the
embedding of flag representations in divided powers is given,
working over an arbitrary finite field $\mathbb{F}$, using the category of functors from
finite-dimensional $\mathbb{F}$-vector spaces to $\mathbb{F}$-vector spaces.
\end{abstract}


\section{Introduction}

Let $V$ be a finite-dimensional $\field$-vector space over a finite
field $\field$; the flag variety of complete flags of length $r$ in $V$
induces a permutation representation $\field [\flag_r] (V)$ of the
general linear group $GL (V)$, which is of interest in representation
theory. The notation is derived from the fact that the flag representation
arises as the evaluation on the space $V$ of a functor $\field
 [\flag_r]$  in  the  category $\f$ of functors from finite-dimensional
 $\field$-vector spaces to $\field$-vector spaces; similarly,  the
 divided power functors $\Gamma ^n$ induce $GL (V)$-representations $\Gamma ^n (V)$. Motivated by
questions from algebraic topology (and  working over the prime field
$\field_2$), Crabb and Hubbuck \cite{ach,ch}
associated to a sequence $\underline{s}$ of integers, $s_1 \geq
\ldots  \geq s_r > s_{r+1} =0$, a
morphism 
\begin{eqnarray}
\field_2 [\flag_r ] (V)
\rightarrow 
\Gamma ^{[\underline{s}]_2} (V)
\end{eqnarray}
or $GL(V)$-modules, 
 where $[\underline{s}]_q$ is the integer $\Sigma _{i=1} ^r (q^{s_i}
 -1)$. The motivating observation of this paper is that this  is
 defined globally as a  natural transformation 
\begin{eqnarray}
\label{eqn:flag-rep}
\phi_{\underline{s}} : 
\field [\flag_r ]
\rightarrow 
\Gamma ^{[\underline{s}]_q}
\end{eqnarray}
in the functor category $\f$ and for each finite field $\field$, where $q = |\field|$.

The Crabb-Hubbuck morphism $\phi_{\underline{s}}$ arises  in the construction of the
ring of lines (developed independently in the dual
situation by Repka and Selick \cite{rs}). Namely,  the divided power
functors form a commutative graded algebra in
$\f$ and  the ring of lines is the graded sub-functor generated by the  images of the morphisms
$\phi_{\underline{s}}$, which  forms a sub-algebra of
$\Gamma ^*$.

The ring of lines  is of
interest in relation to the study of the primitives under
the action of the Steenrod algebra on the singular homology $H_* (BV ;
\field_2)$ of the classifying space of $V$; the primitives arise in a
number of questions in algebraic topology. This relation can be explained from the point of view of the functor category $\f$ as follows, by identifying $H_* (BV; \field_2)$
with the graded vector space $\Gamma ^* (V)$. Steenrod reduced
power operations correspond to natural transformations of the form
$\Gamma ^a \rightarrow \Gamma^b$, $a \geq b$ (see \cite{k_genI}). Define the Steenrod kernel
functors $K^a$ by:
\[
K^a := \mathrm{Ker} \{ \Gamma ^a \rightarrow \bigoplus_{f \in \hom (\Gamma^a,
  \Gamma ^b), a>b} \Gamma ^b \}. 
\]
These functors form a commutative graded algebra in 
$\f$. The primitives in $H_* (BV; \field_2)$ are obtained by evaluating on $V$. The analysis of the primitives  is dual to the 
study of the indecomposables for the action of the Steenrod reduced
powers on $H^* (BV; \field_2)$; this is a difficult problem
which has attracted much interest. For the field
$\field_2$, the complete structure is known only for spaces of
dimension at most four; the case of dimension three is due to Kameko
and a published account is available in Boardman
\cite{kam,board}; Kameko  also announced the case of dimension four,
which  has since been calculated by Sum, a student of Nguyen Hung.

The morphism $\phi_{\underline{s}}$ of equation  (\ref{eqn:flag-rep})  maps to
$K^{[\underline{s}]_2}$ for elementary reasons and  the ring of lines is
a graded sub-algebra in $\f$ of  $K^*$.  A fundamental question is to determine in which degrees the
ring of lines  coincides with the Steenrod kernel. Motivated by this
question, Crabb and Hubbuck \cite[Proposition 3.10]{ch} gave an explicit criterion upon
the sequence $(s_i)$ with respect to the dimension of $V$ for the
morphism $\phi_{\underline{s}}$  to be a monomorphism.

The purpose of this note is two-fold; to present a proof exploiting
the  category $\f$ and to
generalize the result to an arbitrary finite field $\field$,
with $q = |\field |$. The main result of the
paper is the following:

\begin{THM}
 Let $r$ be a natural number and  $\underline{s} = (s_1 > \ldots > s_r> s_{r+1}=0) $ be a sequence of
 integers which satisfies the condition
$[s_i - s_{i+1}]_q\geq  (q-1)(\dim V  -i +1)$, for $1 \leq i \leq  r$; then the
 morphism $\phi_{\underline{s}}$ induces a monomorphism 
\[
 \field [\flag_r](V)
\hookrightarrow 
\Gamma ^{[\underline{s}]_q} (V)
.
\]
\end{THM}

The  proof sheds light upon the method proposed by Crabb and
Hubbuck;  namely, the proof establishes the stronger result that the
composite with a morphism induced by the iterated diagonal on divided power
algebras and the Verschiebung morphism is a
monomorphism. This is of interest in the light of recent work over the
field $\field_2$ by Grant Walker, Reg Wood \cite{ww} and Tran Ngoc Nam
generalizing the result of Crabb and Hubbuck, relaxing the required
hypothesis on the sequence $(s_i)$.

These techniques can be used to provide further information on the nature of the embedding results; for instance:

\begin{PROP}
Let $\underline{s}$ be a sequence of integers $(s_1> \ldots > s_r
>s_{r+1} =0)$ and $V$ be a finite-dimensional vector space for which
the morphism $\phi_{\underline{s}}(V) : \field [\flag_ r ] (V)
\rightarrow \Gamma ^{[\underline{s}]_q} (V) $ is a monomorphism. 

Let $\underline{s^+}$ denote the sequence given by $s_i^+ = s_i +1$, 
for $1 \leq i \leq r$,  and $s^+_{r+1}=0$. Then the morphism 
\[
\phi_{\underline{s^+}}(V) : 
\field [\flag_ r ] (V)
\rightarrow 
\Gamma ^{[\underline{s^+}]_q} (V) 
\]
is a monomorphism. 
\end{PROP}

\bigskip
\noindent
{\bf Acknowledgement :}
 The author would like to thank Grant Walker for conversations which
 led to the current approach to the result of Crabb and Hubbuck and  for
 his interest.


\tableofcontents

\section{Preliminaries}

Fix a finite field  $\field=\field_q$, and write  $q=p^m$, where $p$ is
the characteristic of $\field$. Let $\f$ be the category of functors from finite-dimensional
 $\field$-vector spaces to $\field$-vector spaces. The group of units
$\field^\times$ is isomorphic to the cyclic group $\zed /(q-1)\zed$ via $i
\mapsto \lambda ^i$, for a generator $\lambda$. In particular, the
group has order prime  to $p$, hence the category of $\field[\field
^\times]$-modules is semi-simple. This gives rise to the weight
splitting of $\f$ (as in \cite{k_genI}):
\[
\f 
\cong 
\prod_{i \in \zed/(q-1)\zed}
\f^{i},
\]
where $\f^{i}$ is the full subcategory of functors such that $F
(\lambda 1_V)= \lambda ^i F(1_V)$ for all finite-dimensional spaces $V$ and $\lambda \in
\field ^\times$. By reduction mod $q-1$, the weight category $\f^k$
can be taken to be defined for $k$ an integer.

The duality functor  $D : \f\op \rightarrow \f$ is defined by $DF
 (V):=F(V^*)^*$ and  $D$ restricts to a functor $D : (\f^k)\op \rightarrow \f^k$

\subsection{Divided powers and the Verschiebung}

 Recall that $\Gamma ^k$ denotes
the $k$th divided power functor, defined as the invariants $\Gamma ^k := (T^k ) ^{\sym_k}$
under the action of the symmetric group permuting the factors of
$T^k $, the $k$th tensor power functor. (By convention, the divided power $\Gamma
^0$ functor is the constant functor $\field$ and $\Gamma^i =0$ for $i <0$). The divided power functor
$\Gamma^k$ is dual  to the
$k$th symmetric power functor $S^k$ and the functors $\Gamma ^k$, $T^k$,
$S^k$ all belong to the weight category $\f^{{k}}$. 

The divided power functors $\Gamma ^*$ form a graded exponential
functor; namely for finite-dimensional vector spaces $U$, $V$ and a
natural number $n$, there is a binatural isomorphism
\[
\Gamma ^n (U \oplus V) 
\cong 
\bigoplus_{i+j =n}
\Gamma ^i (U) \otimes \Gamma ^j (V)
.
\]
 This has important consequences  (see \cite{ffss}, for example); in particular,  for pairs of natural numbers $(a,b)$,   there are cocommutative coproduct morphisms
$
 \Gamma^{a+b}
\stackrel{\Delta}{\rightarrow} 
\Gamma ^a \otimes \Gamma ^b
$
 and  commutative product morphisms 
$
 \Gamma ^a \otimes \Gamma ^b 
\stackrel{\mu}{\rightarrow}
 \Gamma ^{a+b},
$
which are  coassociative  (respectively associative) in the
appropriate graded sense.

The Verschiebung is  a natural surjection 
 $
\mathcal{V} :
\Gamma ^{qn}
\twoheadrightarrow 
\Gamma ^n, 
 $ 
for integers $n \geq 0$,  dual to the
Frobenius $q$th power  map on symmetric powers. More
generally there
is a Verschiebung  morphism $\mathcal{V}_p : \Gamma
^{np}\twoheadrightarrow \Gamma^{(1)}$, dual to the Frobenius $p$th power
map, where $(-)^{(1)}$ denotes the Frobenius twist functor 
(see \cite{ffss}). The Verschiebung $\mathcal{V}$ is obtained by iterating
$\mathcal{V}_p$ $m$ times, where $q=p^m$.

The $p$th truncated symmetric power functor $\overline{S}^n$ is
defined by imposing the relation $v^p = 0$; similarly the $q$th
truncated symmetric power functor $\tilde{S}^n$ is given by forming
the quotient by the relation $v^q =0$.  There is a natural surjection
$\tilde{S}^n \twoheadrightarrow \overline{S}^n$;  over a prime field the functors coincide. 

Dualizing gives the following:

\begin{defn}
For $n$ a natural number, 
\begin{enumerate}
\item
let $\tilde{\Gamma}^n$ denote the kernel of the composite
\[
\Gamma ^n 
\stackrel{\Delta}{\rightarrow} 
\Gamma ^{n-q} \otimes \Gamma ^q 
\stackrel{1 \otimes \mathcal{V}}{\rightarrow} 
\Gamma ^{n-q} \otimes \Gamma ^1;
\]
\item
let $\overline{\Gamma}^n$ be the kernel of the composite 
\[
\Gamma ^n 
\stackrel{\Delta}{\rightarrow} 
\Gamma ^{n-p} \otimes \Gamma ^p 
\stackrel{1 \otimes \mathcal{V}_p}{\rightarrow} 
\Gamma ^{n-p} \otimes (\Gamma ^1)^{(1)}. 
\]
\end{enumerate}
\end{defn}

\begin{lem}
\label{lem:Gamma-bar}
Let $n$ be a natural number.
 \begin{enumerate}
\item
The functor $\tilde{\Gamma}^n$ is dual to $\tilde{S}^n$.
\item
The functor $\overline{\Gamma}^n$ is isomorphic to $\overline{S}^n$ and is
     simple.
\item
\label{item:tilde-Gamma-vanish}
$\tilde{\Gamma}^n(\field^d)=0$ if $d \leq (n-1)/(q-1)$. 
\end{enumerate}
\end{lem}

\begin{proof}
The first statement follows from the definitions and the identification $DS^n
 = \Gamma ^n$. The simplicity of $\overline{S}^n$ is a standard fact \cite{k_genI};
 the isomorphism follows from the fact that the simple functors of $\f$
 are self-dual  \cite{k_genII}.

The final statement is an elementary verification, which is a consequence of the observation that the maximal degree of a free
 $q$-truncated symmetric algebra on $d$ variables of degree one is $d(q-1)$.
\end{proof}

\subsection{Further properties of divided powers}

\begin{nota}
For $s$ a natural number, let $[s]_q$ denote the integer $q^s -1$ and, for $\underline{s}$ a sequence of integers, $s_1 \geq \ldots  \geq s_r > s_{r+1} =0$, let $[\underline{s}]_q$ denote $\Sigma _i [s_i]_q$.
\end{nota}

\begin{nota}
The element of $\Gamma ^k (V)$ corresponding to the symmetric tensor
 $x^{\otimes k} \in T^k (V)$, for $x$ an element of $V$, will be denoted simply by $x^{\otimes k}$.
\end{nota}

\begin{lem}
\label{lem:product-expansion}
Let $s$ be a positive integer. 
For any $x \in V$, the class $x^{\otimes [s]_q} \in \Gamma ^{[s]_q} (V)$
is equal to the product 
\[
\prod_{j=0}^{sm-1} 
x^{\otimes (p-1) p^j},
\]
where $|\field|= q = p^m$.
\end{lem}

This Lemma is a consequence of the following well-known general property of the
divided power functors. 

\begin{lem}
\label{lem:gamma-splitting}
Let $a_0, \ldots, a_t$ and $0=r_0 < r_1 < \ldots < r_t$ be sequences of
 natural numbers such that, for $0 \leq i < t$, $a_i p^{r_i} <p^{r_{i+1}}$. Then the
 composite
\[
 \Gamma^{\Sigma_{i=0 }^t a_i p^{r_i}}
\rightarrow 
\bigotimes_{i=0}^t \Gamma ^{a_i p^{r_i}}
\rightarrow
 \Gamma^{\Sigma_{i=0 }^t a_i p^{r_i}}
\]
is an isomorphism, where the first morphism is the coproduct and the
 second the product.
\end{lem}

\begin{proof}
Using the (co)associativity of the product (respectively the coproduct),
 it suffices to prove the result for $t=1$. The composite morphism
 in this case is multiplication by the scalar $\binom{a_0 +
 a_1p^{r_1}}{a_0}$, which is equal to one modulo $p$, since $a_0 < p^{r_1}$, by hypothesis. 
\end{proof}

The following Lemma is the key  to the construction of the Crabb-Hubbuck morphism, $\phi_{\underline{s}}$, in Section \ref{sect:chmorph}.

\begin{lem}
\label{lem:product-vanish}
Let $x$ be an element of $V$, $s$ be a natural number  and $0<i \leq [s]_q$ be an integer. The product 
$x^{\otimes [s]_q} x ^{\otimes i}$ is zero in $\Gamma ^{[s]_q +i} (V)$.
\end{lem}

\begin{proof}
The element $x^{\otimes [s]_q}$ is equal to the product $\prod_{j=0}^{sm-1} 
x^{\otimes (p-1) p^j}$, by Lemma \ref{lem:product-expansion}. Similarly,
 considering the $p$-adic expansion $i= \Sigma_{j=0}^{sm-1}i_j p^j$,
 where $0 \leq i_j < p$ and at least one $i_j$ is non-zero, there is an
 equality $x^{\otimes i}= \prod _{j=0}^{sm-1} 
x^{\otimes i_j p^j}$. The product is associative and commutative, hence
 it suffices to show that, if $i_j\neq 0$, then $x^{\otimes
 (p-1)p^j}x^{\otimes i_j p^j}$ is zero. Up to non-zero scalar in
 $\field_p^\times$, the element $x^{\otimes i_j p^j}$ is the $i_j$-fold
 product of $x^{\otimes p^j}$ (since $0< i_j < p$, by hypothesis), hence it suffices to show that $x^{\otimes
 (p-1)p^j}x^{\otimes p^j}$ is zero. This element identifies with $\binom
{p^{j+1}}{p^j}x ^{\otimes p^{j+1}}$ and the scalar is zero in $\field$.

\end{proof}

The behaviour of the Verschiebung morphism with respect to products is
important.

\begin{lem}
\label{lem:versch-product}
Let $\beta_1, \ldots , \beta_k$ be positive integers such that
 $\Sigma_{i=1}^k \beta_i = p N$ for some integer $N$. The composite 
\[
 \bigotimes_{i=1}^k \Gamma ^{\beta_i}
\stackrel{\mu}{\rightarrow}
\Gamma ^{pN}
\stackrel{\mathcal{V}_p}{\rightarrow }
(\Gamma ^N)^{(1)},
\]
where $\mathcal{V}_p$ is the Verschiebung and $\mu$ is the product,
 is trivial unless for each $i$, $\beta_i =
 p \beta'_i$, $\beta'_i\in \mathbb{N}$. In this case, there is a commutative diagram
\[
 \xymatrix{
 \bigotimes_{i=1}^k \Gamma ^{p \beta'_i}
\ar[r]^{\mu}
\ar[d]_{\bigotimes \mathcal{V}_p}
&
\Gamma ^{pN}
\ar[d]^{\mathcal{V}_p}
\\
  \bigotimes_{i=1}^k (\Gamma ^{\beta'_i})^{(1)}
\ar[r]_{\mu^{(1)}}
&
(\Gamma ^N)^{(1)}.
}
\]
\end{lem}

\begin{proof}
The statement is more familiar in the dual situation, where it
 corresponds to the fact that the diagonal of the symmetric power
 algebra commutes with the Frobenius.
\end{proof}

\begin{rem}
An analogous statement holds for iterates of $\mathcal{V}_p$ and for  the Verschiebung $\mathcal{V}:
 \Gamma ^{qN}\rightarrow \Gamma^N$.
\end{rem}


\section{Projectives and flags}

This section introduces the flag  functors and relates them to the standard projective generators of the category $\f$.

\subsection{The projective and flag functors}
For $r$ a natural number,  the standard projective object
$P_{\field^r} $ in $\f$ is given by $P_{\field^r} (V)= \field [\hom
(\field ^r, V)]$ and is determined up
to isomorphism by $\hom_\f (P_{\field^r}, G ) = G (\field^r)$. In
particular, the projective $P_\field$ is the functor $V \mapsto \field [V]$.

 The weight
splitting determines a direct sum decomposition 
\[
P_\field \cong \bigoplus _{i \in \zed/ (q-1) \zed} P_\field ^i
\]
in which $P_\field ^i$  is indecomposable for $i \neq 0$ and $P_\field
^0$ admits a decomposition $P_\field ^0 = \field \oplus
\overline{P_\field^0}$ (Cf \cite[Lemma 5.3]{k_genI}, noting that Kuhn
uses a splitting associated to the multiplicative semigroup $\field$).

There is a Künneth isomorphism for projectives so that, for $r$ a
positive integer, $P_\field ^{\otimes r}$ is projective and identifies
with the projective functor.

\begin{defn}
\label{def:flags}
For $r$ a positive integer, let  $\field [\flag _r]$ be the functor  which is defined
in terms of complete flags of length $r$ as follows. 

As a vector space, $\field [\flag_r](V)$ has basis the set of complete flags
of length $r$. A morphism $V\rightarrow W$ sends a complete flag to its  image, if this is a complete flag, and  to zero otherwise.  
\end{defn}

Recall that a functor $F$ of $\f$ is said to be constant-free if $F(0)=0$.

\begin{lem}
\label{lem:flag-properties}
Let $r\geq s >0 $ be  integers.
\begin{enumerate}
\item
 The functor $\field [\flag_r]$ is
constant-free and belongs to $\f ^0$. 
\item
\label{item:flag-diag}
There is a diagonal morphism in $\f$:
\[
 \field [\flag_r ]
\rightarrow 
\field [\flag _r ]\otimes 
\field [\flag_r].
\]
\item
\label{item:flag-quotient}
There is a surjection
\[
\pi_{r,s} : \field [\flag _r ]
\twoheadrightarrow 
\field [\flag_{s}]
\]
which forgets the subspaces of dimension greater than $s$.
\end{enumerate}
\end{lem}

\subsection{Structure of the projectives}

\begin{lem}
\label{lem:unique-morphism}
Let $n \geq 1$ be an integer, then $\hom_\f (\overline{P_\field^0}, 
\Gamma ^{n (q-1)}) = \field$.
\end{lem}

\begin{proof}
The result follows from the Yoneda lemma, the weight splitting and the fact that
$n \geq 1$ allows passage to the constant-free part, $\overline{P_\field^0}$.
\end{proof}

Recall that a functor is said to be finite if it has a finite
composition series.

\begin{prop}
\label{prop:overlineP}
\cite{k_genI,k_genII} 
\begin{enumerate}
\item
The surjection $P_\field \twoheadrightarrow \field [\flag_1]$ induces
an isomorphism $\overline{P_\field ^0} \cong \field [\flag_1]$. 
\item
There is an inverse system of finite functors $ \ldots \rightarrow \qpol_k
\overline{P_\field^0} \rightarrow \qpol_{k-1}\overline{P_\field^0}
\rightarrow \ldots $ 
such that 
\begin{enumerate}
\item
$
\overline{P_\field^0} 
\cong
\lim_\leftarrow
\qpol_k
\overline{P_\field^0}
$;
\item
for $k\geq 1$, 
$ \qpol_k
\overline{P_\field^0}$
is isomorphic to the image of any non-trivial morphism $\overline{P_\field^0}
\rightarrow \Gamma ^{k (q-1)}$;
\item
for $k \geq 2$, 
there is a non-split short exact sequence 
\[
0 \rightarrow \tilde{\Gamma}^{k(q-1)}
\rightarrow
\qpol_k
\overline{P_\field^0}
\rightarrow 
\qpol_{k-1}
\overline{P_\field^0}
\rightarrow 
0.
\]
\end{enumerate}
\item
The functor $\overline{P_\field ^0}$ is dual to a locally-finite functor.
\item
If $\field$ is the prime field $\field_p$, then 
the functor  $\overline{P_\field ^0}$ is uniserial with composition
factors $\{ \overline{\Gamma}^{k(p-1)} | k \geq 1\}$, each occurring
with multiplicity one. 
\end{enumerate}
\end{prop}

\begin{proof}
It is more straightforward to deduce the result from the description of
 the dual $D\overline{P_\field^0}$; this is isomorphic
 to the functor $(\bigoplus_{k=0}^\infty S^{k(q-1)})/ \langle v^q -v
 \rangle$, where the relation is induced by the weight zero part of the
 ideal $\langle v^q - v \rangle$ (this is deduced from  \cite[Lemma
 4.12]{k_genI} by applying the evident weight splitting).
\end{proof}

It is important to have a measure of how good an approximation  $\qpol_k \overline{P_\field^0}$ is to
$\overline{P_\field^0}$.

\begin{lem}
\label{lem:factorize}
Let $k \geq 1$ be an integer. Up to scalar in $\field^\times$, there is a unique non-trivial morphism
 $\overline{P_\field^0}\rightarrow \Gamma ^{k(q-1)}$. Any non-trivial
 morphism $\overline{P_\field^0}\rightarrow \Gamma ^{k(q-1)}$
\begin{enumerate}
\item
factors as 
\[
 \overline{P_\field^0}
\twoheadrightarrow 
\qpol_k \overline{P_\field^0}
\hookrightarrow 
\Gamma ^{k(q-1)}
\]
and 
\item
\label{item:iso-bound}
induces a  monomorphism 
\[
 \overline{P_\field^0} (V)\hookrightarrow \Gamma ^{k(q-1)} (V)
\]
if $\dim V \leq k$.
\end{enumerate}
\end{lem}

\begin{proof}
The unicity follows from Lemma \ref{lem:unique-morphism} and
 the factorization follows from this unicity together with the
 identification of $\qpol_k \overline{P_\field^0}$ which is given in
 Proposition \ref{prop:overlineP}.

The kernel of the  surjection $\overline{P_\field^0}\twoheadrightarrow \qpol_k
 \overline{P_\field^0}$ has a filtration with  subquotients of
 the form $\tilde{\Gamma}^{l(q-1)}$ with $l > k $. The functor
 $\tilde{\Gamma}^{l(q-1)}$ is zero when evaluated on spaces with $\dim V \leq k$, by Lemma \ref{lem:Gamma-bar} (\ref{item:tilde-Gamma-vanish}). It
 follows that the kernel is zero  when evaluated on such spaces. 
 Thus $\overline{P_\field^0}(V)\rightarrow \qpol_k
 \overline{P_\field^0}(V)$ is an isomorphism
 when $\dim V \leq k$ and  the
 result follows from the first part of the Lemma.
\end{proof}


\section{The flag morphism of Crabb and Hubbuck}
\label{sect:chmorph}

Throughout this section, let $r$ denote a fixed positive integer and  $\underline{s}$ denote a fixed decreasing sequence of
positive integers, $s_1 \geq s_2 \geq \ldots \geq s_r >s_{r+1}=0$.

\begin{nota}
For each positive integer $s$, let  $\phi_s$ be the element 
of  $\hom _\f (\overline{P_\field^0}, \Gamma ^{[s]_q} )$ which sends the
canonical generator of $\overline{P_\field^0} (\field) \cong \field
[\flag_1](\field)$ to $\iota^{\otimes [s]_q}$, where $\iota$ is any
generator of $\field$. (The morphism is independent of the
choice of $\iota$).
\end{nota}

\begin{defn}
For $\underline{s}$ a sequence of positive integers, let
$\tilde{\phi}_{\underline{s}}$ denote the morphism 
\[
\tilde{\phi} _{\underline{s}}:
P_{\field^r} \cong  
P_\field ^{\otimes r}
\stackrel{\otimes \phi_{s_i}}{\longrightarrow }
\bigotimes_{i=1}^r \Gamma ^{[s_i]_q}
\rightarrow 
\Gamma ^{[\underline{s}]_q}
\]
in which the second morphism is induced by the product.
\end{defn}

The following Proposition is proved  in the case $q=2$ in \cite{ch}.

\begin{prop}
\label{prop:ch-morphism}
The morphism $\tilde{\phi} _{\underline{s}}$ factorizes as 
\[
  P_\field ^{\otimes r}
\twoheadrightarrow
\field[\flag _r ]
\stackrel{\phi_{\underline{s}}}{\rightarrow}
\Gamma ^{ [\underline{s}]_q}.
\]
\end{prop}

\begin{proof}
Fixing a basis of $\field ^r$, a canonical basis element  of $P_{\field^r} (V)$ is an
ordered sequence $(v_i)$ of $r$ elements of $V$. The morphism
$\tilde{\phi}_{\underline{s}}$ sends this generator to $\prod_{i=1}
^{r} v_i^{\otimes [s_i]_q} = \prod_{i=1}^r \prod_{j=0}^{s_i-1} v_i
^{\otimes (q-1)q^j}$.

Using this notation, define a natural surjection $P_{\field^r}
\twoheadrightarrow \field [\flag_r] $ by sending $(v_i)$ to the flag
$\langle v_1\rangle < \langle v_1, v_2\rangle<  \ldots  <  \langle v_1, \ldots , v_r\rangle$ if the
elements are linearly independent and zero otherwise. The proposition
asserts that $\tilde{\phi}_{\underline{s}}$ factorizes across this surjection.
 
The result follows as in the proof of \cite[Lemma 3.1]{ch}, by applying 
Lemma \ref{lem:product-vanish}.
\end{proof}


\section{The embedding theorem}

The purpose of this section is to prove the main result of the paper,
stated here as  Theorem
\ref{thm:embed-flag},  which gives a criterion  for 
\[
\phi_{\underline{s}} (V)
:
\field [\flag_r](V)
\rightarrow 
\Gamma ^{[\underline{s}]_q} (V)
\]
to be a monomorphism, where
 $r$ is a positive integer and $\underline{s}= s_1 > s_2 > \ldots > s_r >s_{r+1}= 0$ is  a strictly decreasing sequence
 of integers.

\begin{rem}
 The morphism $\phi_{\underline{s}}$ is clearly not a monomorphism of functors,
since the functor $\Gamma ^{[\underline{s}]_q}$ is finite  whereas $\field
[\flag_r ]$ is highly infinite. 
\end{rem}

However, Lemma \ref{lem:factorize}   (\ref{item:iso-bound}) provides the key calculational input, which is restated as the following:

\begin{lem}
\label{lem:case-r1} 
Let $s$ be a natural number. The morphism
 $\phi_s : \field [\flag_1] \cong \overline{P_\field ^0} \rightarrow \Gamma
 ^{[s]_q}$ induces  a monomorphism $\phi_s (V)$ if $[s]_q\geq (q-1)\dim V$. 
\end{lem}

The theorem is proved by an induction using Lemma \ref{lem:case-r1} to provide the inductive step. 
 The strategy involves composing $\phi_{\underline{s}}$ with a morphism $\delta_{\underline{s}}
$ (defined below) to give a morphism $\psi_{\underline{s}}$ which is amenable to induction. The key to setting up the induction is Lemma  \ref{lem:key-step}.

Write $[\underline{s}]_q = r [s_r]_q + \Sigma _{i=1} ^{r-1} ([s_i]_q - [s_{r}]_q )$; thus 
the coproduct gives  a morphism 
 $
\Delta : 
\Gamma ^{[\underline{s}]_q}
\rightarrow 
\Gamma ^{r [s_r]_q}
\otimes 
\Gamma ^{\Sigma ([s_i]_q - [s_r]_q)}. 
 $ 
For each $i$, $[s_i]_q- [s_r]_q= q^{s_r}[s_i - s_r]_q$, 
 hence  there is an iterated Verschiebung morphism 
 $
\mathcal{V}^{s_r} :  \Gamma ^{\Sigma ([s_i]_q - [s_r]_q)}
\twoheadrightarrow
\Gamma ^{\Sigma [s_i -s_r]_q}.
$

\begin{defn}
\label{def:delta-gamma-theta}
\ 
\begin{enumerate}
\item
\label{def:delta}
Let 
$
 \delta_{\underline{s}}
 :
\Gamma ^{[\underline{s}]_q}
\rightarrow 
\Gamma ^{r [s_r]_q}
\otimes 
\Gamma ^{\Sigma_{i=1}^{r-1}[s_i -s_r]_q} 
$
be the composite morphism
\[
\Gamma ^{[\underline{s}]_q}
\stackrel{\Delta}{\rightarrow}
\Gamma ^{r [s_r]_q}
\otimes 
\Gamma ^{\Sigma ([s_i]_q - [s_r]_q)}
\stackrel{1 \otimes \mathcal{V}^{s_r}}{\longrightarrow}
\Gamma ^{r [s_r]_q}
\otimes 
\Gamma ^{\Sigma_{i=1}^{r-1}[s_i -s_r]_q}.
\]
\item
\label{def:psi}
Let 
$
 \psi_{\underline{s}} :
\field [\flag_r ]
\rightarrow
\Gamma ^{r [s_r]_q}
\otimes 
\Gamma ^{\Sigma_{j=1}^{r-1}[s_j -s_r]_q} 
 $
be the composite morphism $\delta_{\underline{s}}\circ\phi _{\underline{s}} $.
\end{enumerate}
\end{defn}

The following elementary observation is recorded as a Lemma.

\begin{lem}
\label{lem:phi-psi}
Let $V$ be a finite-dimensional vector space. If the morphism
$\psi_{\underline{s}}(V)$ is a monomorphism, then
$\phi_{\underline{s}} (V)$ is a monomorphism.
\end{lem}

\begin{nota}
Let $\underline{s'}$ denote the sequence (of length $r-1$) of positive integers 
$(s_1 - s_r > \ldots > s_{r-1}- s_r >0)$.
\end{nota}

The following Lemma is  the key to the inductive proof, and relies upon the
fact that the iterated Verschiebung is used in the definition of $\psi
 _{\underline{s}}$. Observe that the Crabb-Hubbuck
morphism associated to the sequence of integers $(s_r, \ldots ,
s_r)$ of length $r$ induces a morphism 
\[
\phi_{(s_r, \ldots, s_r)} : \field [\flag_r  ]
\rightarrow 
\Gamma ^{r[s_r]_q}.
\] 

\begin{lem}
\label{lem:key-step}
The morphism $\psi_{\underline{s}}$ identifies with the composite
 morphism 
\[
\xymatrix{
\field [\flag _r ]
\ar[r]^(.3){\mathrm{diag}}
&
\field [\flag _r ]\otimes\field [\flag_r]
\ar[rr]^{1 \otimes \pi_{r,r-1}}
&
\ 
&
\field [\flag _r ]\otimes\field [\flag_{r-1}]
\ar[d]^{\phi _{(s_r, \ldots, s_r)} \otimes \phi_{\underline{s'}}}
\\
&
&
&
\Gamma ^{r[s_r]_q}
\otimes 
\Gamma ^{\Sigma_{j=1}^{r-1}[s'_j]_q} 
}
\]
\end{lem}

\begin{proof}
We are required to prove that the following diagram is commutative.
\[
 \xymatrix{
\field [\flag _r ]
\ar[r]^(.3){\mathrm{diag}}
\ar[d]_{\phi_{\underline{s}}}
&
\field [\flag _r ]\otimes\field [\flag_r]
\ar[rr]^{1 \otimes \pi_{r,r-1}}
&&
\field [\flag _r ]\otimes\field [\flag_{r-1}]
\ar[d]^{\phi _{(s_r, \ldots, s_r)} \otimes \phi _{\underline{s'}}}
\\
\Gamma ^{[\underline{s}]_q}
\ar[rrr]_{\delta_{\underline{s}}}
&&
&
\Gamma ^{r[s_r]_q}
\otimes 
\Gamma ^{\Sigma_{j=1}^{r-1}  [s'_j]_q}.
}
\]
Choose a  surjection $P_{\field^r}\twoheadrightarrow \field [\flag_r]$ as in
the proof of Proposition \ref{prop:ch-morphism}; it is equivalent to
prove the commutativity of the diagram  obtained by composition, replacing the top left entry by $P_{\field
 ^r}$. The analysis of the composite morphism 
\begin{eqnarray}
\label{eqn:composite}
 P_{\field^r}
\rightarrow 
 \Gamma ^{[\underline{s}]_q}
\rightarrow
\Gamma ^{r[s_r]_q}
\otimes 
\Gamma ^{\Sigma_{j=1}^{r-1}  [s'_j]_q}.
\end{eqnarray}
around the bottom of the diagram can then  be carried out as follows.  

The definition of the morphism $\delta_{\underline{s}}$ and of the morphism $\phi_{\underline{s}}$ implies that this composite factors across  
\[
\bigotimes_{i=1}^r 
\Gamma^{[s_i]_q}
\stackrel{\mu}{\rightarrow}
\Gamma^{[\underline{s}]_q}
\stackrel{\Delta}{\rightarrow}
\Gamma ^{r[s_r]_q} \otimes \Gamma ^{q^{s_r}\Sigma _{j=1}^{r-1}[s'_j]_q}
\]
where $\mu$ denotes the product on divided powers and $\Delta$ the diagonal. 

The exponential algebra structure of $\Gamma ^*$ (essentially the fact that these functors take values in bicommutative Hopf algebras) implies that there is a  commutative diagram 
\[
 \xymatrix{
\bigotimes _{i=1}^r P_{\field}
\ar@/^2pc/@<2ex>[rr]^{\tilde{\phi}_{\underline{s}}}
\ar[r]^{\bigotimes _{i=1}^r \phi_{s_i}}
&
\bigotimes_{i=1}^r 
\Gamma^{[s_i]_q}
\ar[r]^\mu
\ar[d]_{\bigoplus \Delta}
&
\Gamma^{[\underline{s}]_q}
\ar[d]^\Delta
\ar@/^3pc/@<5ex>[dd]^{\delta_{\underline{s}}}
\\
&
\bigoplus 
\Big \{
\bigotimes _{i=1}^r 
\big(
\Gamma^{\alpha_i} \otimes \Gamma ^{\beta_i}
\big)
\Big \}
\ar[r]^{\bigoplus (\mu \otimes \mu)}
\ar@{-->}[rd]_{\xi}
&
\Gamma ^{r[s_r]_q} \otimes \Gamma ^{q^{s_r}\Sigma _{j=1}^{r-1}[s'_j]_q}\ 
\ar[d]^{1 \otimes \mathcal{V}^{s_r}}
\\
&
&
\Gamma ^{r[s_r]_q} \otimes \Gamma ^{\Sigma _{j=1}^{r-1}[s'_j]_q}
}
\]
where the sum is labelled over
 sequences of pairs of natural numbers $(\alpha_i, \beta_i)$ satisfying 
 $\alpha_i + \beta _i = [s_i]_q$ for each $i$ and $\Sigma \alpha_i = r[s_r]_q$.

Consider the composite morphism $\xi$  in the diagram; Lemma  \ref{lem:versch-product}  implies that the only
 components of this morphism which are non-trivial are those
 corresponding to sequences $(\alpha_i, \beta_i)$ for which $\beta_i =
 q^{s_r}\beta'_i$ for natural numbers $\beta_i'$. The condition
 $\alpha_i + q^{s_r}\beta'_i= [s_i]_q$ implies that $\alpha_i$ is
 non-zero; it follows that  $\alpha_i \geq [s_r]_q$, for each $i$, since $\alpha_i \equiv [s_r]_q \mod (q^{s_r})$. The
 condition $\Sigma_i \alpha_i = r[s_r]_q$ therefore implies that
 $\alpha_i = [s_r]_q$ for each $i$. It follows that $\xi$ has only one
 non-zero component and a straightforward verification shows that the
composite  corresponds to the composite around the top of the diagram
given in the statement of the Lemma.
\end{proof}

The inductive argument is simplified using the following:

\begin{lem}
\label{lem:inductive-step}
Let $V$ be a finite-dimensional vector space for which the morphism
$\phi_{\underline{s}'}(V)$ is a monomorphism. Then the morphism
$\psi_{\underline{s}}( V)$ is a monomorphism if and only if the
composite morphism 
\[
\xymatrix{
\field [\flag _r ]
\ar[r]^(.3){\mathrm{diag}}
&
\field [\flag _r ]\otimes\field [\flag_r]
\ar[rr]^{\phi_{s_r, \ldots, s_r}\otimes \pi_{r,r-1}}
&
\ 
&
\Gamma ^{r[s_r]_q}
\otimes
\field [\flag_{r-1}]
}
\]
induces a monomorphism when evaluated upon $V$.
\end{lem}

\begin{proof}
This  follows from the identification of $\psi_{\underline{s}}$
which is given in Lemma \ref{lem:key-step}.
\end{proof}

Using the fact that $\field[\flag_{r-1}](V) $ is generated by complete flags of length $r-1$, this allows the decomposition into components.

\begin{nota}
For $\fl$ a complete flag  of length $r-1$ in $V$, let 
\begin{enumerate}
\item
$
\langle 
\fl
\rangle
\leq V
$ denote the $(r-1)$-dimensional subspace of $V$ defined by $\fl$;
\item
$\field [\flag_r]_\fl(V)$ denote the subspace generated by flags containing $\Phi$;
\item
$\gamma_{\fl}$ denote the image of $[\fl] \in \field[\flag_{r-1}]
(V)$ under the Crabb-Hubbuck morphism $\phi_{s_r, \ldots , s_r} (V) :
\field [\flag_{r-1}](V)\rightarrow \Gamma ^{(r-1) [s_r]_q} (V)$.
\end{enumerate}
\end{nota} 

\begin{rem}
 The space $\field [\flag_r]_\fl(V)$ is isomorphic to $\field [\flag_1] (V/ \langle \fl \rangle)$.
\end{rem}

\begin{lem}
\label{lem:fl-reduction}
Let $V$ be a finite-dimensional vector space for which the morphism
$\phi_{\underline{s}'}(V)$ is a monomorphism. The morphism
$\psi_{\underline{s}} (V)$ is a monomorphism if and only if, for each
complete flag $\fl$ in $V$ of length $r-1$, the 
restriction of 
\[
\field [\flag_r] (V) 
\stackrel{\phi_{s_r, \ldots, s_r}}{\longrightarrow}
\Gamma ^{r[s_r]_q} (V)
\]
to $
\field [\flag_r]_\fl (V) 
$
is a monomorphism.
\end{lem}

\begin{proof}
A straightforward consequence of Lemma \ref{lem:inductive-step}.
\end{proof}

\begin{nota}
Let $V$ be a finite-dimensional vector space and  $\fl$ be a complete
flag in $V$ of length $r-1$. Let $\rho _\fl$ denote the  composite linear map
\[
\field [\flag_r]_\fl (V) 
\stackrel{\cong}{\rightarrow}
\field [\flag_1](V/ \langle \fl \rangle)
\stackrel{\phi_{s_r}}{\rightarrow} 
\Gamma ^{[s_r]_q} (V/  \langle \fl \rangle)
\]
induced by  the projection $V \twoheadrightarrow V/
\langle \fl \rangle$ and the morphism $\phi_{s_r}$. 
\end{nota}

\begin{nota}
For $V$  a finite-dimensional vector space,  $\fl$  a complete
flag in $V$ of length $r-1$ and $\sigma$  a section of the
projection $V \twoheadrightarrow V / \langle \fl \rangle$, let
\[
\gamma_\fl \cap _{\sigma} 
:\Gamma ^{[s_r]_q} (V/ \langle \fl \rangle
) 
\rightarrow \Gamma ^{r[s_r]_q} (V) 
\]
denote the linear morphism  induced by the section $\sigma$
  followed by the product with $\gamma_\fl$ with respect to the
  algebra structure of $\Gamma ^* (V)$.
\end{nota}

\begin{lem}
\label{lem:mono}
Let $V$, $\fl$ and $\sigma$ be as above.  The linear
morphism
\[
\gamma_\fl \cap _{\sigma} 
:
\Gamma ^{[s_r]_q} (V/ \langle \fl \rangle) 
\rightarrow 
\Gamma ^{r[s_r]_q} (V)
\]
is a monomorphism. 
\end{lem}

\begin{proof}
The result follows from  the exponential structure of the divided
power functors, since the element $\gamma_\fl$ is the image of an
element of $\Gamma ^{(r-1)[s_r]_q } (\langle \fl \rangle )$ under the
morphism induced by the natural inclusion.
\end{proof}

\begin{lem}
\label{lem:composite}
Let $V$, $\fl$, $\sigma$ be as above. The
restriction of $\phi_{s_r, \ldots, s_r}(V)$ to $\field [\flag_r ]_\fl
(V)$ identifies with  the linear morphism $(\gamma_\fl \cap_\sigma) \circ  
\rho_{\fl}$.
\end{lem}

\begin{proof}
The result follows from the definition of the  morphism $\phi_{s_r, \ldots, s_r} (V)$.
\end{proof}

Lemmas \ref{lem:mono} and \ref{lem:composite} together imply the following result:

\begin{lem}
\label{lem:monomorphism-equivalent}
Let $V$ be a finite-dimensional vector space and   $\fl$ be a complete
flag in $V$ of length $r-1$. The restriction of $\phi_{s_r, \ldots
  ,s_r} (V)$ to $\field [\flag_r]_\fl (V)$ is a monomorphism if and
only if 
\[
\phi_{s_r} :\field [\flag_1](V/ \langle \fl \rangle)
\rightarrow 
\Gamma ^{[s_r]_q} (V/  \langle \fl \rangle)
\]
is a monomorphism. 
\end{lem}

\begin{rem}
By lemma \ref{lem:case-r1},  a sufficient condition is
\[
[s_r]_q
\geq (q-1) \dim (V/ \langle \fl \rangle) = (q-1) (\dim V -r +1) .
\] 
When $q=2$, this is an equivalent condition.
\end{rem}

Putting these results together, one obtains the following generalization of \cite[Proposition 3.10]{ch}.

\begin{thm}
\label{thm:embed-flag}
Suppose that the sequence $\underline{s} $ satisfies the condition
$[s_i - s_{i+1}]_q\geq  (q-1)(\dim V  -i +1)$, for $1 \leq i \leq  r$. Then the
 morphism $\phi_{\underline{s}}$ induces a monomorphism 
\[
 \field [\flag_r](V)
\hookrightarrow 
\Gamma ^{[\underline{s}]_q} (V)
.
\]
\end{thm}

\begin{proof}
The result is proved by induction upon $r$, starting with the initial
 case, $r=1$, which is provided by Lemma \ref{lem:case-r1}. For the
 inductive step, by Lemma \ref{lem:phi-psi}, it is sufficient to show
 that $\psi_{\underline{s}}(V)$ is a monomorphism, under the given
 hypotheses.

 Observe that the hypotheses upon $\underline{s}$ imply
 that $\underline{s}'$ also satisfy the hypotheses with respect to $V$,
 so that the morphism $\phi_{\underline{s}'}(V)$ is injective, by induction.
 Hence Lemma \ref{lem:fl-reduction} reduces the proof to showing
 that the restriction of $\phi_{s_r, \ldots, s_r}$ to $
\field [\flag_r]_\fl (V) 
$ is a monomorphism, for each complete flag $\fl$ of length $r-1$ in
 $V$. The inductive step is completed by combining Lemma
 \ref{lem:monomorphism-equivalent} with Lemma \ref{lem:case-r1}.
\end{proof}


\section{A stabilization result}

The techniques employed in the proof of Theorem \ref{thm:embed-flag}
can be used to provide further information on the nature of the
embedding results. For instance, one has a direct  proof of the
following stabilization result. 

\begin{prop}
Let $\underline{s}$ be a sequence of integers $(s_1> \ldots > s_r
>s_{r+1} =0)$ and $V$ be a finite-dimensional vector space for which
the morphism $\phi_{\underline{s}}(V) : \field [\flag_ r ] (V)
\rightarrow \Gamma ^{[\underline{s}]_q} (V) $ is a monomorphism. 

Let $\underline{s^+}$ denote the sequence given by $s_i^+ = s_i +1$, 
for $1 \leq i \leq r$,  and $s^+_{r+1}=0$. Then the morphism 
\[
\phi_{\underline{s^+}}(V) : 
\field [\flag_ r ] (V)
\rightarrow 
\Gamma ^{[\underline{s^+}]_q} (V) 
\]
is a monomorphism. 
\end{prop}

\begin{proof}
The diagonal induces a morphism $\Gamma ^{[\underline{s^+}]_q}
\rightarrow \Gamma ^{ q [\underline{s}]_q} \otimes \Gamma ^{(q-1)
  r}$. Hence, composing with the Verschiebung on the first morphism
gives $\eta : \Gamma ^{[\underline{s^+}]_q}
\rightarrow \Gamma ^{[\underline{s}]_q} \otimes \Gamma ^{(q-1)
  r}$, as in Definition \ref{def:delta-gamma-theta}.

There is a commutative diagram 
\[
\xymatrix{
\field [\flag_r] 
\ar[rr]^{\phi_{\underline{s^+}}}
\ar[d]_{\mathrm{diag}}
&&
 \Gamma ^{[\underline{s^+}]_q}
\ar[d]^\eta 
\\
\field [\flag_r] \otimes \field [\flag _r] 
\ar[rr]_{\phi_{\underline{s}}\otimes \phi_{1, \ldots , 1}}
&&
\Gamma ^{[\underline{s}]_q} \otimes \Gamma ^{(q-1) r},
}
\]
the  commutativity of which  is established by an argument similar
to that employed in the proof of Lemma \ref{lem:key-step}.

It suffices to show that the composite 
\[
\field [\flag_r] (V) 
\stackrel{\phi_{\underline{s^+}}} {\longrightarrow }
\Gamma ^{[\underline{s^+}]_q}(V)
\stackrel{\eta}{\rightarrow }
\Gamma ^{[\underline{s}]_q} \otimes \Gamma ^{(q-1) r}(V)
\]
is a monomorphism. By hypothesis, the morphism $\phi _{\underline{s}} (V)$ is a
monomorphism. As in the inductive step of the proof of Theorem
\ref{thm:embed-flag}, the result then follows from the fact that the
morphism $\field [\flag_r ] (V) \rightarrow \Gamma ^{(q-1)r} (V)$ is
non-trivial. The latter follows from the fact that the hypothesis upon
$\phi_{\underline{s}} (V)$ implies that $V$ has dimension at least
$r$. 
\end{proof}

\begin{rem}
This argument is related to standard
techniques using the Kameko $Sq^0$ operation \cite{kam}, which is
based in an essential way upon the analysis of the Verschiebung
morphism. 
\end{rem}



\def\cftil#1{\ifmmode\setbox7\hbox{$\accent"5E#1$}\else
  \setbox7\hbox{\accent"5E#1}\penalty 10000\relax\fi\raise 1\ht7
  \hbox{\lower1.15ex\hbox to 1\wd7{\hss\accent"7E\hss}}\penalty 10000
  \hskip-1\wd7\penalty 10000\box7}
\providecommand{\bysame}{\leavevmode\hbox to3em{\hrulefill}\thinspace}
\providecommand{\MR}{\relax\ifhmode\unskip\space\fi MR }
\providecommand{\MRhref}[2]{%
  \href{http://www.ams.org/mathscinet-getitem?mr=#1}{#2}
}
\providecommand{\href}[2]{#2}

\end{document}